\documentclass[]{article}
\usepackage{amssymb,enumerate}
\usepackage{amsmath}
\usepackage{graphics}
\usepackage[utf8]{inputenc}

\newtheorem{theorem}{Theorem}[section]
\newtheorem{corollary}[theorem]{Corollary}

\newtheorem{remark}[theorem]{Remark}
\newtheorem{proposition}[theorem]{Proposition}

\title{Heun functions related to entropies}
\author{Adina  B\u{a}rar\thanks{Technical University of Cluj-Napoca, Department of Mathematics, Memorandumului Street 28, 400114, Cluj-Napoca, Romania}, Gabriela Raluca Mocanu\thanks{Romanian Academy, Cluj-Napoca Branch, Astronomical Institute, Cire\c{s}ilor Street 19, 400487, Cluj-Napoca, Romania,
gabriela.mocanu$@$academia-cj.ro}, Ioan Ra\c{s}a\thanks{Technical University of Cluj-Napoca, Department of Mathematics \newline, Memorandumului Street 28, 400114, Cluj-Napoca, Romania}}

\date{}
\begin{document}

\maketitle

Subject Class: 33E30, 33C05, 94A17
\newline \newline
\textbf{Abstract}
\newline

We consider the indices of coincidence for the binomial, Poisson, and negative binomial distributions. They are related in a simple manner to the R\'{e}nyi entropy and Tsallis entropy. We investigate some families of Heun functions containing these indices of coincidence. For the involved Heun functions we obtain closed forms, explicit expressions, or representations in terms of hypergeometric functions.

\section{Introduction}
Consider the general Heun equation
\begin{equation}
u''(x) + \left ( \frac{\gamma}{x} + \frac{\delta}{x-1}+\frac{\epsilon}{x-a} \right )u'(x) + \frac{\alpha \beta x - q}{x(x-1)(x-a)} u(x) = 0,\label{eq:1.1}
\end{equation}
where $a \notin \{ 0, 1\}$, $\gamma \notin \{ 0,-1,-2, \dots \}$ and $\alpha + \beta + 1 = \gamma + \delta + \epsilon$. The solution $u(x)$ normalized by $u(0)=1$ is called the \emph{(local) Heun function} and is denoted by $Hl(a,q;\alpha , \beta; \gamma , \delta ; x)$. For details see, e.g., \cite{twos}, \cite{2} and the references therein.

Consider also the confluent Heun equation
\begin{equation}
u''(x) + \left ( 4p+\frac{\gamma}{x} + \frac{\delta}{x-1} \right )u'(x) + \frac{4 p\alpha x - \sigma}{x(x-1)} u(x) = 0,\label{eq:1.2}
\end{equation}
where $p \neq 0$. Its solution normalized by $u(0)=1$ is called the \emph{confluent Heun function} and is denoted by $HC(p,\gamma, \delta, \alpha, \sigma;x)$.

For a positive integer $n$ consider the binomial probability distribution
\begin{equation*}
{n \choose k} x^k (1-x)^{n-k}, \quad x \in [0,1], \quad k=0,1,\dots , n,
\end{equation*}
the negative binomial distribution
\begin{equation*}
{n+k-1 \choose k} x^k (1+x)^{-n-k}, \quad x \in [0,+\infty ), \quad k=0,1,\dots ,
\end{equation*}
and the Poisson distribution
\begin{equation*}
e^{-nx} \frac{(nx)^k}{k!}, \quad x \in [0,+\infty ), \quad k=0,1,\dots \quad.
\end{equation*}

The corresponding indices of coincidence were denoted in~\cite{5} by
\begin{equation}
F_n (x): = \sum _{k=0}^n \left ( {n \choose k} x^k (1-x)^{n-k}\right )^2, \label{eq:1.3}
\end{equation}
\begin{equation}
G_n(x): = \sum _{k=0}^\infty \left ( {n+k-1 \choose k} x^k (1+x)^{-n-k} \right )^2, \label{eq:1.4}
\end{equation}
\begin{equation}
K_n(x):=\sum _{k=0}^\infty \left ( e^{-nx}\frac{(nx)^k}{k!} \right )^2 . \label{eq:1.5}
\end{equation}

If $S_n(x)$ is an index of coincidence, then $R_n(x) : = - \log S_n(x)$ and $T_n(x):= 1-S_n(x)$ are the R\'{e}nyi, respectively the Tsallis entropies of order $2$.

It was proved in~\cite[(34),(60),(68)]{5} that
\begin{equation}
F_n (x)= Hl \left ( \frac{1}{2}, -n;-2n,1;1,1;x \right ), \label{eq:1.6}
\end{equation}
\begin{equation}
G_n(-x) = Hl \left ( \frac{1}{2}, n;2n,1;1,1;x \right ), \label{eq:1.7}
\end{equation}
\begin{equation}
K_n(x)=HC \left ( n,1,0,\frac{1}{2},2n;x \right ). \label{eq:1.8}
\end{equation}

See also~\cite[(5),(6),(21)]{1}.

In this paper we investigate some families of Heun functions containing $F_n(x)$, $G_n(-x)$ and $K_n(x)$ as particular cases. From the general results we deduce properties of these three functions. Clearly, each property of $F_n(x)$, $G_n(x)$ and $K_n(x)$ can be translated as a property of the associated R\'{e}nyi and Tsallis entropies.

In Section~\ref{sect:2} we use the results of~\cite{3} in order to study the derivatives of some Heun functions. Our main results are concerned with the closed forms of the functions $Hl \left ( \frac{1}{2}, -2n\theta; -2n, 2\theta; \gamma,\gamma;x \right )$ and $Hl \left ( \frac{1}{2}, 2n\theta; 2n, 2\theta; \gamma,\gamma;x \right )$ which generalize the functions $F_n(x)$, respectively $G_n(-x)$.

Section~\ref{sect:3} is devoted to the representation of $Hl \left ( \frac{1}{2}, q; 2q, 1; 1,1;x \right )$ in terms of hypergeometric functions, in the spirit of~\cite{9}.

In Section~\ref{sect:4} we obtain explicit expressions of the function \\ $HC \left (n, j+1, 0, j+\frac{1}{2}, 2n(2j+1);x \right )$ for integers $n\geq 1$ and $j\geq 0$.

Some of the results were also presented in~\cite{1}.

We need the following important formula (see line 3 in Table 2 of \cite{4}):
\begin{eqnarray}
&&Hl(a,q;\alpha, \beta; \gamma, \delta; x) = \label{eq:1.9}\\&& \left ( 1-\frac{x}{a} \right )^{-\alpha - \beta + \gamma + \delta} Hl(a,q-\gamma (\alpha + \beta - \gamma - \delta); -\alpha + \gamma + \delta , -\beta + \gamma + \delta ; \gamma, \delta; x).\nonumber
\end{eqnarray}

\section{Derivatives of Heun functions\label{sect:2}}

From~\eqref{eq:1.1} we get immediately $u'(0) = \frac{q}{a\gamma}$, i.e.,
\begin{equation}
\frac{\rm{d}Hl(a,q;\alpha, \beta;\gamma,\delta;x)}{\rm{dx}}(0) = \frac{q}{a\gamma}.\label{eq:2.1}
\end{equation}

The derivatives of Heun functions were investigated in~\cite{3}. Here we are interested in the case when
\begin{equation}
q=a\alpha \beta.\label{eq:2.2}
\end{equation}

Combining\eqref{eq:2.1} with (14)-(17) from~\cite{3} we get the following result.
\begin{proposition}\label{prop:2.1}
With the above notation,
\begin{equation}
\frac{d}{dx}Hl(a,a\alpha\beta;\alpha,\beta;\gamma,\delta;x) = \frac{\alpha\beta}{\gamma}\left ( 1-\frac{x}{a} \right )Hl(a,q_1;\alpha _1, \beta _1;\gamma +1, \delta +1;x),\label{eq:2.3}
\end{equation}
where
\begin{eqnarray*}
\alpha _1 \beta _1 &=& \alpha \beta + 2 (\gamma + \delta + \epsilon +1),\\
\alpha _1 + \beta _1 &=& \gamma + \delta + \epsilon +3,\\
q_1 &=& a (\alpha \beta + \gamma + \delta) + \gamma + \epsilon +1.
\end{eqnarray*}

\begin{equation}
\frac{d}{dx}Hl(a,a\alpha\beta;\alpha,\beta;\gamma,\delta;x) = \frac{\alpha\beta}{\gamma}\left ( 1-\frac{x}{a} \right )^{-\epsilon} Hl(a,q_2;\alpha _2, \beta _2;\gamma +1, \delta +1;x),\label{eq:2.4}
\end{equation}
where
\begin{eqnarray*}
\alpha _2 \beta _2 &=& \alpha \beta + (\gamma + \delta ) (1- \epsilon),\\
\alpha _2 + \beta _2 &=& \gamma + \delta - \epsilon +1,\\
q_2 &=& a (\alpha \beta + \gamma + \delta) - \gamma  \epsilon .
\end{eqnarray*}
\end{proposition}

\begin{remark}
The equality of the right-hand side members of \eqref{eq:2.3} and \eqref{eq:2.4} is also a consequence of \eqref{eq:1.9}.
\end{remark}

From Proposition~\ref{prop:2.1} we get immediately

\begin{corollary}\label{cor:2.3}
Let $a=1/2$ and $\delta = \gamma$. Then

\begin{eqnarray}
&&\frac{d}{dx}Hl \left ( \frac{1}{2}, \frac{1}{2}\alpha \beta; \alpha , \beta ; \gamma , \gamma ; x \right )\label{eq:2.5} \\ & =& \frac{\alpha \beta}{\gamma} (1-2x)Hl \left ( \frac{1}{2}, \frac{1}{2}(\alpha +2)(\beta +2);\alpha +2, \beta +2; \gamma +1 , \gamma +1; x \right ), \nonumber
\end{eqnarray}

\begin{eqnarray}
&&\frac{d}{dx}Hl \left ( \frac{1}{2}, \frac{1}{2}\alpha \beta; \alpha , \beta ; \gamma , \gamma ; x \right )\label{eq:2.6} \\ & =& \frac{\alpha \beta}{\gamma} (1-2x)^{2\gamma -\alpha -\beta -1}Hl \left ( \frac{1}{2}, \frac{1}{2}(2\gamma-\alpha)(2\gamma - \beta );2\gamma - \alpha , 2\gamma -\beta ; \gamma +1 , \gamma +1; x \right ). \nonumber
\end{eqnarray}
\end{corollary}

\begin{remark}\label{rem:2.4}
The Heun functions in the above right-hand sides satisfy the condition of type \eqref{eq:2.2}, so that it is possible to express their derivatives in terms of other Heun functions.
\end{remark}

Let $k\geq 0$ be an integer and $r$ be a real number. We shall use the notation
\begin{equation*}
(r)_k:=r(r+1)\dots (r+k-1), \quad \mbox{if } k>0; \quad (r)_0:=1.
\end{equation*}

\begin{theorem}\label{thm:2.5}
Let $n\geq 0$ be an integer and $\theta \in \mathbb{R}$. Then
\begin{equation}\label{eq:2.7}
Hl \left ( \frac{1}{2},-2n\theta;-2n,2\theta;\gamma,\gamma;x \right ) = \sum _{k=0}^n 4^k {n \choose k} \frac{(\theta)_k}{(\gamma)_k}(x^2-x)^k .
\end{equation}
\end{theorem}

\subsection*{Proof}
We shall prove \eqref{eq:2.7} by induction with respect to $n$. By using \eqref{eq:1.1} it is easy to see that
\begin{equation*}
Hl \left ( \frac{1}{2}, 0; 0, 2\theta;\gamma, \gamma; x\right )=1,
\end{equation*}
and so \eqref{eq:2.7} is verified for $n=0$. Suppose that it is verified for a given $n\geq 0$. Then, according to \eqref{eq:2.5} we get
\begin{eqnarray*}
\frac{d}{dx}Hl \left ( \frac{1}{2}, -2(n+1)\theta; -2(n+1),2\theta;\gamma,\gamma;x \right ) = \\= \frac{4(n+1)\theta}{\gamma}(2x-1)Hl\left ( \frac{1}{2}, -2n(\theta+1); -2n,2(\theta+1);\gamma+1,\gamma+1;x \right )=\\
= \frac{4(n+1)\theta}{\gamma}(2x-1)\sum _{k=0}^n 4^k {n \choose k} \frac{(\theta +1)_k}{(\gamma +1)_k}(x^2-x)^k =\\
= \sum _{k=0}^n 4^{k+1} (k+1) {n+1 \choose k+1} \frac{(\theta)_{k+1}}{(\gamma)_{k+1}}(x^2-x)'(x^2-x)^k.
\end{eqnarray*}

Taking into account the normalization at $0$, it follows that
\begin{eqnarray*}
Hl \left ( \frac{1}{2}, -2(n+1)\theta; -2(n+1),2\theta;\gamma,\gamma;x \right ) = \\= 1+\sum _{k=0}^n 4^{k+1} {n+1 \choose k+1} \frac{(\theta)_{k+1}}{(\gamma)_{k+1}}(x^2-x)^{k+1} = \\
=\sum _{k=0}^{n+1} 4^k {n+1 \choose k} \frac{(\theta )_k}{(\gamma )_k}(x^2-x)^k,
\end{eqnarray*}
and so \eqref{eq:2.7} is valid for $n+1$. This concludes the proof.

\subsection*{Example}
In particular, for $\theta = \gamma$ \eqref{eq:2.7} yields
\begin{equation*}
Hl \left ( \frac{1}{2},-2n\gamma; -2n, 2\gamma;\gamma, \gamma; x\right )= (1-2x)^{2n}.
\end{equation*}

\begin{corollary}\label{cor:2.7}
Let $\gamma$ and $n$ be integers, $0<\gamma \leq n$, and $\theta \in \mathbb{R}$. Then
\begin{eqnarray}
Hl \left ( \frac{1}{2}, 2n\theta; 2n,2\theta;\gamma,\gamma;x \right ) =\label{eq:2.8}\\= (1-2x)^{-2(n-\gamma+\theta)}\sum _{k=0}^{n-\gamma} 4^k {n-\gamma \choose k} \frac{(\gamma-\theta )_k}{(\gamma )_k}(x^2-x)^k. \nonumber
\end{eqnarray}
\end{corollary}

\subsection*{Proof}
From \eqref{eq:1.9} and \eqref{eq:2.7} we get successively
\begin{eqnarray*}
Hl \left ( \frac{1}{2}, 2n\theta; 2n,2\theta;\gamma,\gamma;x \right ) =\\= (1-2x)^{-2(n-\gamma+\theta)}Hl \left ( \frac{1}{2}, 2(\gamma-n)(\gamma-\theta); -2(n-\gamma),2(\gamma-\theta);\gamma,\gamma;x \right )=\\
= (1-2x)^{-2(n-\gamma+\theta)} \sum _{k=0}^{n-\gamma} 4^k {n-\gamma \choose k} \frac{(\gamma-\theta )_k}{(\gamma )_k}(x^2-x)^k,
\end{eqnarray*}
and this proves \eqref{eq:2.8}.

\subsection*{Example}
For $\gamma = n$, \eqref{eq:2.8} produces
\begin{eqnarray*}
Hl \left ( \frac{1}{2}, 2n\theta; 2n,2\theta;n,n;x \right ) =(1-2x)^{-2\theta}.
\end{eqnarray*}

\begin{corollary}
For $n\geq 1$ we have
\begin{equation}\label{eq:2.9}
F_n(x) = \sum _{k=0}^{n}{n \choose k} {2k \choose k}(x^2-x)^k,
\end{equation}
\begin{equation}\label{eq:2.10}
G_n(x) = (1+2x)^{1-2n}\sum _{k=0}^{n-1} {n-1 \choose k}{2k \choose k}(x^2+x)^k.
\end{equation}
\end{corollary}

\subsection*{Proof}
To prove~\eqref{eq:2.9} it suffices to apply \eqref{eq:1.6} and \eqref{eq:2.7} with $\gamma=1$, $\theta =1/2$. \eqref{eq:2.10} is a consequence of \eqref{eq:1.7} and \eqref{eq:2.8}.

From \eqref{eq:2.9} it is easy to obtain
\begin{equation}\label{eq:2.11}
F_n(x) = \sum _{j=0}^n (1-2x)^{2j} 4^{-j} {n \choose j} \sum _{i=0}^{n-j} \left ( -\frac{1}{4} \right )^i {n-j \choose i} {2i+2j \choose i+j}.
\end{equation}

Similarly, \eqref{eq:2.10} leads to
\begin{equation}\label{eq:2.12}
G_n(x) = \sum _{j=0}^{n-1} (1+2x)^{2j-2n+1} 4^{-j} {n-1 \choose j} \sum _{i=0}^{n-j-1} \left ( -\frac{1}{4} \right )^i {n-j-1 \choose i} {2i+2j \choose i+j}.
\end{equation}

(\eqref{eq:2.12} is also a consequence of \eqref{eq:2.11} and \cite[Cor. 13]{5}).

The following formula was established in \cite{6}:
\begin{equation}\label{eq:2.13}
F_n(x) = \sum _{j=0}^n (1-2x)^{2j} 4^{-n} {2j \choose j}  {2n-2j \choose n-j}.
\end{equation}

Moreover, according to \cite[(56)]{5},
\begin{equation}\label{eq:2.14}
G_n(x) = \sum _{j=0}^{n-1} (1+2x)^{2j-2n+1} 4^{1-n} {2n-2j-2 \choose n-j-1} {2j \choose j}.
\end{equation}

Starting from \eqref{eq:2.13} we get easily
\begin{equation}\label{eq:2.15}
F_n(x) = \sum _{k=0}^n (x^2-x)^k 4^{k-n}  \sum _{j=k}^{n}{j \choose k}{2j \choose j}{2n-2j \choose n-j}.
\end{equation}

Comparing \eqref{eq:2.9} with \eqref{eq:2.15}, and then \eqref{eq:2.11} with \eqref{eq:2.13} we obtain the following (seemingly nontrivial) combinatorial identities:
\begin{equation*}
\sum _{j=k}^n{j \choose k}{2j \choose j}{2n-2j \choose n-j} = 4^{n-k}{n \choose k}{2k \choose k}, \quad 0\leq k \leq n,
\end{equation*}

\begin{equation*}
\sum _{i=0}^{n-j}\left ( -\frac{1}{4} \right )^i {n-j \choose i} {2i+2j \choose i+j} = 4^{j-n}{2j \choose j}{2n-2j \choose n-j}{n \choose j}^{-1}, \quad 0\leq j \leq n.
\end{equation*}

\begin{remark}\label{rem:2.10}
Equations \eqref{eq:2.9}-\eqref{eq:2.15} provide explicit expressions for the indices of coincidence $F_n(x)$ and $G_n(x)$, which are also the Heun functions \\ $Hl(1/2, -n; -2n,1;1,1;x)$, respectively $Hl(1/2,n;2n,1;1,1;-x)$. As mentioned in the Introduction, they can be used in order to study the corresponding R\'{e}nyi and Tsallis entropies; for details see also \cite{1}, \cite{5}, \cite{7}, \cite{8}, where bounds and shape properties for $F_n(x)$, $G_n(x)$ and the entropies are obtained.
\end{remark}

\section{Heun functions and hypergeometric functions\label{sect:3}}

Let us consider the Heun function $Hl \left ( \frac{1}{2}, q; 2q,1;1,1;x \right )$, which was studied also in the preceding section. The parameters $\alpha = 2q$, $\beta = 1$, $\gamma = 1$, $\delta = 1$, $\epsilon = 2q$ satisfy the condition (33) from \cite{9}, namely
\begin{equation*}
a=\frac{1}{2}, \quad \gamma + \delta =2, \quad q = a\alpha \beta + a (1-\delta)\epsilon.
\end{equation*}

Consequently, according to \cite[(35) and (36)]{9},
\begin{equation}
Hl \left ( \frac{1}{2},q;2q,1;1,1;x \right ) = \frac{u(x)}{u(0)},\label{eq:3.1}
\end{equation}
where
\begin{equation}
u(x) = \sum _{k=0}^\infty \frac{(1/2)_k(q)_k}{k! (q+1/2)_k} \frac{q}{q+k} {_2F_1}(2q,1;1+2q+2k;x),\label{eq:3.2}
\end{equation}
and
\begin{equation*}
u(0) = {_3F_2} \left ( \frac{1}{2},q,q;q+\frac{1}{2},q+1;1 \right ).
\end{equation*}

Here $q\notin \{ 0, -1, \dots \}$ $\cup$ $\left \{ -\frac{1}{2}, -\frac{3}{2}, \dots \right \}$, ${_2F_1}$ is the Gauss hypergeometric function and ${_3F_2}$ is the Clausen hypergeometric function (see \cite{9}).

Let $q=m/2$, where $m\geq 1$ is an integer. We shall determine a closed form of the hypergeometric function ${_2F_1} (m,1;m+2k+1;x)$ which appears in \eqref{eq:3.2}.

First,

\begin{equation*}
{_2F_1}(m,1;m+2k+1;x) = 1 + \sum _{j=1}^\infty \frac{(m)_j}{(m+2k+1)_j}x^j,
\end{equation*}
so that
\begin{equation*}
\left ( x^{m+2k}{_2F_1}(m,1;m+2k+1;x) \right )^{(2k+1)} = (m)_{2k+1}\frac{x^{m-1}}{1-x}.
\end{equation*}

Now it is a matter of calculus to find
\begin{eqnarray}
&&{_2F_1}(m,1;m+2k+1;x) = \label{eq:3.3} \\ &=& (m)_{2k+1}x^{-m-2k} \left ( \frac{(1-x)^{2k}}{(2k)!} \left ( e_{2k}-\log (1-x) \right ) - \sum_{j=0}^{2k}a_{jk}x^j - \sum _{i=0}^{m-2} \frac{x^{i+2k+1}}{(i+1)_{2k+1}} \right ),\nonumber
\end{eqnarray}
where
$e_n:=1+\frac{1}{2}+ \dots + \frac{1}{n}$ and
\begin{equation*}
a_{jk}:=\frac{1}{(2k)!} \left ( \sum_{i=0}^{j-1} {2k \choose i} \frac{(-1)^i}{j-i} + (-1)^j {2k \choose j} e_{2k}\right ).
\end{equation*}

In particular,
\begin{equation}
{_2F_1}(1,1;2;x) = -\frac{\log (1-x)}{x}\label{eq:3.4}
\end{equation}
and
\begin{equation}
{_2F_1}(2,1;3;x) = - 2 \frac{x+\log (1-x)}{x^2}.\label{eq:3.5}
\end{equation}

\section{Confluent Heun functions\label{sect:4}}

In this section we obtain explicit expressions for some confluent Heun functions.

Let $\gamma \neq 0$. From \eqref{eq:1.2} we get $u'(0) = -\sigma / \gamma$. This fact, combined with \cite[(21)]{10}, leads to
\begin{equation}
\frac{d}{dx}HC(p,\gamma, 0, \alpha, 4p\alpha;x) = -\frac{\sigma}{\gamma} HC (p, \gamma +1, 0, \alpha +1, 4p(\alpha +1);x), \label{eq:4.1}
\end{equation}
\begin{equation}
\frac{d}{dx}HC(p,\gamma, 0, \alpha, 4p\alpha;x) = \frac{\sigma}{\gamma}(x-1)HC(p,\gamma +1, 2, \alpha +2, 4p (\alpha +1)-\gamma -1; x). \label{eq:4.2}
\end{equation}

According to \eqref{eq:1.8}, the index of coincidence $K_n(x)$ is the confluent Heun function $HC\left (n,1,0,1/2,2n;x\right )$.

From \eqref{eq:1.8} and \eqref{eq:4.2} we obtain
\begin{equation}
HC \left ( n, 2, 2, \frac{5}{2}, 6n-2;x \right ) = \frac{1}{2n(x-1)}K_n'(x).\label{eq:4.3}
\end{equation}

The main result of this section is contained in
\begin{theorem}\label{thm:4.1}
Let $n\geq 1$ and $j\geq 0$ be integers. Then
\begin{equation}
HC\left ( n,j+1,0,j+\frac{1}{2},2n(2j+1);x \right ) = (-n)^{-j} {2j \choose j}^{-1}K_n^{(j)}(x).\label{eq:4.4}
\end{equation}
\end{theorem}

\subsection*{Proof}
As a confluent Heun function, $K_n(x)$ satisfies the following equation (see also \cite[(68)]{5}):
\begin{equation*}
xK''_n(x) + (4nx+1)K'_n(x)+2nK_n(x)=0.
\end{equation*}

Taking the j$^{th}$ derivative we get
\begin{equation*}
\left ( K_n^{(j)} \right )''(x) + \left (4n+\frac{j+1}{x} \right ) \left (  K_n^{(j)}\right )'(x) + \frac{2n(2j+1)}{x}K_n^{(j)}(x)=0.
\end{equation*}

According to \eqref{eq:1.2}, this is a confluent Heun equation, and so
\begin{equation}
\frac{K_n^{(j)}(x)}{K_n^{(j)}(0)} = HC\left ( n,j+1,0,j+\frac{1}{2},2n(2j+1);x \right ). \label{eq:4.5}
\end{equation}

Let us remark that \eqref{eq:4.5} can be proved alternatively by using \eqref{eq:1.8}, \eqref{eq:4.1} and induction after $j$.

It remains to find the value of $K_n^{(j)}(0)$. In \cite[Th. 4]{6} it was proved that
\begin{equation*}
K_n^{(j)}(x) = \frac{2}{\pi} 4^j (-n)^j \int _0 ^{\frac{\pi}{2}} e^{-4nx(\sin {t})^{2j}}dt.
\end{equation*}

It follows that
\begin{equation*}
K_n^{(j)}(0) = \frac{2}{\pi}4^j(-n)^j \int _0^{\frac{\pi}{2}} (\sin t)^{2j} dt = \frac{2}{\pi} 4^j (-n)^j \pi 2^{-2j-1}{2j\choose j},
\end{equation*}
i.e.,
\begin{equation}
K_n^{(j)}(0) = (-n)^j {2j\choose j}.\label{eq:4.6}
\end{equation}

Now \eqref{eq:4.4} is a consequence of \eqref{eq:4.5} and \eqref{eq:4.6}; this concludes the proof.

\section{Final remarks and future work}

a) It is easy to see that in \eqref{eq:2.3} one can take $\alpha_1 = \alpha+2$ and $\beta_1 = \beta+2$. Hence \eqref{eq:2.3} can be written as
\begin{eqnarray}
&&\frac{d}{dx}Hl(a,a\alpha\beta;\alpha, \beta, \gamma, \delta; x) = \nonumber \\&& \frac{\alpha \beta}{\gamma} \left (1-\frac{x}{a} \right ) Hl(a,q_1;\alpha+2,\beta+2;\gamma+1,\delta +1;x). \label{eq:5.1}
\end{eqnarray}

Clearly \eqref{eq:5.1} is similar to the identity
\begin{equation*}
\frac{d}{dx}{_2F_1} (a,b;c;x) = \frac{ab}{c}{_2F_1}(a+1,b+1;c+1;x)
\end{equation*}
used in order to increment the parameters of a Gauss hypergeometric function. Other such identities are more complicated, for example

\begin{eqnarray}
(1-x)^{1-a}\frac{d^m}{dx^m} \left [ (1-x)^{a+m-1} {_2F_1} (a,b;c;x) \right ] \nonumber\\
= \frac{(-1)^m(a)_m(c-b)_m}{(c)_m} {_2F_1}(a+m,b;c+m;x).\label{eq:5.2}
\end{eqnarray}

It would be interesting to explore the existence of similar identities that increment/decrement, individually, the parameters $\alpha$, $\beta$, $\gamma$, $\delta$ of the Heun function $Hl$.

b) \eqref{eq:3.3} and other similar formulas can be obtained if one combines identities like \eqref{eq:5.2} and \eqref{eq:3.4}. It is a pleasant calculation to derive \eqref{eq:3.5} in this manner.

c) From Remark \ref{rem:2.4} it follows that \eqref{eq:2.5} and \eqref{eq:2.6} can be iterated and can be combined. This possibility will be exploited in a forthcoming paper. Here we offer only a sample result (see also \cite[(16)]{1}:

\begin{equation*}
Hl \left ( \frac{1}{2}, (i-n)(2i+1);2(i-n), 2i+1;i+1,i+1;x \right )
\end{equation*}
\begin{equation*}
= \frac{(2i)!!}{(2i-1)!!}4^{-n}{n \choose i}^{-1} \sum _{j=0}^{n-i} 4^j {i+j \choose i} {2i+2j \choose i+j} {2n-2i-2j \choose n-i-j} \left ( x-\frac{1}{2} \right )^{2j},
\end{equation*}
with $i=0,1, \dots ,n$.

d) In \eqref{eq:4.1} and \eqref{eq:4.2}, $\delta=0$. Therefore, the results of \cite[Section 2]{ones} can be applied. Moreover, equating the right sides of \eqref{eq:4.1} and \eqref{eq:4.2}, we get an identity, similar to \eqref{eq:1.9}, for confluent Heun functions. Such identities deserve a separate study.

\subsection*{Acknowledgements}
GM is partially supported by a grant of the Romanian Ministry of National Education
and Scientific Research, RDI Programme for Space Technology and Advanced Research -
STAR, project number 513, 118/14.11.2016.


\end{document}